\documentclass{amsart}
\usepackage{amssymb,amsmath}
\usepackage{hyperref}

\newcommand{\integers}{\ensuremath{\mathbb{Z}}}

\newcommand{\reals}{\ensuremath{\mathbb{R}}}

 \newcommand{\union}{\cup}

\theoremstyle{plain}
\newtheorem{theorem}{Theorem}[section]

\newtheorem{lemma}[theorem]{Lemma}
\newtheorem{corollary}[theorem]{Corollary}

\newtheorem{defn}{Definition}[section]

\newtheorem{example}{Example}[section]

\title{Boundaries of Hyperbolic Metric Spaces}

\author{Corran Webster \and Adam Winchester}

\address{Department of Mathematical Sciences,
University of Nevada Las Vegas,
Las Vegas, NV 89154}

\email{cwebster@unlv.nevada.edu}

\date{September 18, 2003}



\begin{document}

\begin{abstract}
    We investigate the relationship between the metric boundary and
    the Gromov boundary of a hyperbolic metric space.  We show that
    the Gromov boundary is a quotient topological space of the metric
    boundary, and that therefore a word-hyperbolic group has an
    amenable action on the metric boundary of its Cayley graph.  This
    result has significance for the study of Lip-norms on group
    C*-algebras.
\end{abstract}

\maketitle

\section{Introduction}

The Gromov boundary of a hyperbolic metric space has been extensively
studied, but the Gromov boundary is not guaranteed to exist for
non-hyperbolic metric spaces.  Gromov~\cite{gromov:1} introduced
another boundary which makes sense for any metric space, but this was
little studied until Rieffel~\cite{rieffel:groupalgebras} showed that
this second boundary, called the metric boundary in his papers, is
important in the study of metrics on the state spaces of group
C*-algebras.

If $G$ is a countable discrete group equipped with a length function
$\ell$, and $C^{*}_{r}(G)$ is its reduced C*-algebra, then one has a
seminorm $L_{\ell}(f) = \|[M_{\ell}, f]\|$ defined on a dense
*-subalgebra of $C^{*}_{r}(G)$, where $M_{\ell}$ is multiplication by
$\ell$ and $f$ operates by convolution on $\ell^{2}(G)$.  This in turn
gives a metric on the state space of $C^{*}_{r}(G)$ by
\[
  \rho_{L_{\ell}}(\varphi, \psi) = \sup \{|\varphi(f) - \psi(f)| :
  L_{\ell}(f) \le 1 \},
\]
and a natural question to ask is whether the topology generated by
this metric coincides with the weak-* topology on the state space,
ie.~the seminorm is a
Lip-norm~\cite{rieffel:groupmetrics,rieffel:noncommmetrics,rieffel:compactquantummetrics}.
 Rieffel proves that this is in fact the case for $\integers^{d}$ with
certain length functions, and a critical requirement in his proof is
that the action of $\integers^{d}$ on its metric boundary is always
amenable.

There is some interest, then, in knowing when the action of a group is
amenable on its metric boundary.  In the case of word-hyperbolic
groups with the standard word-length metric, it is known that the
action of a word-hyperbolic group on its Gromov boundary is
amenable~\cite{anantharaman:amenablegroupoids,anantharaman:amenability},
and as Rieffel points out in~\cite{rieffel:groupalgebras}, if there is
an equivariant, continuous surjection from the metric boundary onto
the Gromov boundary, then the action of the group on the metric
boundary must be amenable.

We show that this is in fact the case, and more: the Gromov boundary
is a quotient topological space of the metric boundary in a completely
natural way, and that the quotient map is therefore such an
equivariant, continuous surjection from the metric boundary to the
Gromov boundary.

We note here that Ozawa and Rieffel~\cite{ozawarieffel:hyperbolic}
have shown that, for hyperbolic groups, $L_{\ell}$ is in fact a
Lip-norm using techniques which do not use the notion of the metric
boundary.  However these methods do not work for $\integers^{d}$, and
we hope that our result may lead to a unified way of showing that the
seminorms for these groups are in fact Lip-norms.

This paper is part of an undergraduate research project between the
authors.  The authors would like to thank Michelle Schultz for
organizing the undergraduate research seminar at UNLV, and Marc
Rieffel for encouraging this line of research.

\section{The Gromov Boundary}

There are many different but equivalent definitions for a hyperbolic
metric space, but for our purposes we are only interested in a couple. 
We follow Alonso, Smith, et.~al.~\cite{alonso:wordhyperbolicgroups}, 
in our presentation, and a more complete discussion of hyperbolic 
spaces can be found there.

\begin{defn}
    A metric space $(X,d)$ is \emph{geodesic} if given any two points
    $x$, $y \in X$, there is an isometry $\gamma$ from the
    interval $[0, d(x,y)]$ into $X$.
    
    If $(X,d)$ is a metric space, with a base-point $0$, we define an
    \emph{inner product} by
    \[
      (x \cdot y)_{0} = \frac{1}{2}(d(x,0) + d(y,0) - d(x,y)).
    \]
    Where the base point is implicit, we will just write $(x \cdot
    y)$.
    
    The metric space $(X,d)$ is \emph{hyperbolic} if it is geodesic
    and there is some $\delta \ge 0$ such that
    \begin{equation}\label{eqn:hyperboliccond}
      (x \cdot y) \ge \min \{(x \cdot z), (y \cdot z)\} - \delta
    \end{equation}
    for all $z \in X$.
\end{defn}

One can show that although the constant $\delta$ may be different for
different base-points, whether or not the space is hyperbolic does not
depend on the choice of base-point.

We have a particular interest in groups whose Cayley graphs are
hyperbolic, and there is an equivalent definition based on properties
of generators and relations alone.  We note that if $G$ is a group
with a finite presentation $\langle S | R \rangle$, then given a
reduced word $w$ in the generators, $S$, with $w = e$ in $G$, we can
write $w$ as a product
\[
  w = \prod_{k=1}^{n} u_{k}^{-1}r_{k}u_{k},
\]
where $u_{k}$ is a word in $S \union S^{-1}$, and $r_{k} \in R \union
R^{-1}$.  For a given $w$, let $n_{w}$ be the smallest possible number
of terms in such a product, and let $l(w)$ be the length of $w$.

\begin{defn}
    Let $G$ be a group with a finite presentation $\langle S | R
    \rangle$.  We say that $G$ is \emph{word-hyperbolic} if it
    satisfies a \emph{linear isoperimetric inequality}: there is some
    $K \ge 0$ such that
    \[
      n_{w} \le K l(w),
    \]
    for all reduced words $w$ with $w = e$ in $G$.
\end{defn}

One can show that the choice of generators and relations does not
affect whether or not the group is word-hyperbolic and, moreover, a
group is word-hyperbolic if and only if its Cayley graph (regarded as
a 1-complex with the graph metric) is hyperbolic.

Perhaps the simplest way to consider the Gromov boundary is as the
limit points of geodesic rays, where two geodesic rays are considered
equivalent if they are a finite distance apart.  This definition
highlights similarities between the Gromov boundary and the metric
boundary discussed in the next section.  However, the most useful
definition of the Gromov boundary for our purposes is in terms of the
inner product.

\begin{defn}
    Let $(X,d)$ be a metric space.  We say that a sequence $x_{k}$ 
    \emph{converges to infinity (in the Gromov sense)} if
    \[
      \lim_{n,k \to \infty} (x_{n} \cdot x_{k}) = \infty.
    \]
    Given two sequences $x = (x_{n})_{n=1}^{\infty}$ and $y =
    (y_{n})_{n = 1}^{\infty}$ which both converge to infinity, we
    define a relation $\sim$ by
    \[
      x \sim y \iff \lim_{n \to \infty} (x_{n} \cdot y_{n}) = \infty.
    \]
\end{defn}

If $(X,d)$ is a hyperbolic metric space, then $\sim$ is in fact an
equivalence relation on sequences which converge to infinity.  It is 
worthwhile noting that if $(X,d)$ is hyperbolic then
\[
  x \sim y \iff \lim_{n,k \to \infty} (x_{n} \cdot y_{k}) = \infty.
\]

If $(X,d)$ is not hyperbolic, the relation $\sim$ will not be an 
equivalence relation, in general:

\begin{example}
    Consider the Cayley graph of $\integers^{2}$ with the standard
    generators and relations.  Let $x_{n} = (n,0)$, $y_{n} = (0,n)$
    and $z_{n} = (n,n)$.  All three sequences converge to infinity,
    but although $x \sim z$ and $y \sim z$, $x \not\sim y$.
\end{example}

We define the Gromov boundary $\partial_{G}X$ of a hyperbolic metric
space $(X,d)$ to be the set of equivalence classes of sequences which
converge to infinity.  We will say that a sequence in $X$ converges to
an equivalence class in $\partial_{G}X$ if it is an element of the
equivalence class.

We can topologise the boundary by extending the inner product to
$\overline{X}^{G} = X \union \partial_{G}X$.

\begin{defn}
    Let $(X,d)$ be a hyperbolic metric space, and let $x$, $y \in
    \overline{X}^{G}$.  Then we define
    \[
      (x \cdot y) = \inf \{ \liminf_{n} (x_{n} \cdot y_{n}) : x_{n} \to 
      x,\ y_{n} \to y,\ \text{and}\ x_{n}, y_{n} \in X \}.
    \]
\end{defn}

One can show that if this inner product is restricted to $X$, it is
the same as the original inner product on $X$.  Indeed, if $\omega 
\in \partial_{G}X$, and $y \in X$, we have
\[
  (\omega \cdot y) = \inf \{ \liminf_{n} (x_{n} \cdot y) : x_{n} \to
  \omega,\ \text{and}\ x_{n} \in X \}.
\]
It is also the case that if $(X,d)$ is hyperbolic, with
\[
  (x \cdot y) \ge  \min \{ (x \cdot z), (y \cdot z) \} - \delta
\]
for all $x$, $y$ and $z \in X$, then the same identity holds for this
extended inner product. We have
\[
  (x \cdot y) \ge \min \{ (x \cdot z), (y \cdot z) \} - \delta,
\]
for all $x$, $y$ and $z \in \overline{X}^{G}$.

We then can say that a sequence $x_{n} \in \overline{X}^{G}$ converges
to $\omega \in \partial_{G}X$ if and only if
\[
  (\omega \cdot x_{n}) \to \infty.
\]
With this definition, it can be shown that $\overline{X}^{G}$ is a
compactification of $X$.

\section{The Metric Boundary}

We now consider the metric compactification and the metric boundary. 
The most succinct definition is that the metric compactification
$\overline{X}^{d}$ of a metric space $(X,d)$ corresponds to the pure
states of the commutative, unital, C*-algebra $\mathcal{G}(X,d)$
generated by the functions which vanish at infinity on $X$, the
constant functions, and the functions of the form
\[
  \varphi_{y}(x) = d(x,0) - d(x,y),
\]
where $0$ is some fixed base-point (which does not affect the
resulting algebra).  The metric boundary $\partial_{d}X$ is simply
$\overline{X}^{d} \setminus X$.

More concretely, we can understand the metric boundary as a limit of
rays in much the same way as the simple definition of the Gromov
boundary.

\begin{defn}\label{defn:geodesicrays}
    Let $(X,d)$ be a metric space, and $T$ an unbounded subset of
    $\reals^{+}$ containing $0$, and let $\gamma: T \to X$.  We say
    that
    \begin{enumerate}
	\item $\gamma$ is a \emph{geodesic ray} if
	\[
	  d(\gamma(s), \gamma(t)) = |s - t|
	\]
	for all $s$, $t \in T$.
	
	\item $\gamma$ is an \emph{almost-geodesic ray} if for every
	$\varepsilon > 0$, there is an integer $N$ such that
	\[
	  |d(\gamma(t), \gamma(s)) + d(\gamma(s), \gamma(0)) - t| <
	  \varepsilon
	\]
	for all $t$, $s \in T$ with $t \ge s \ge N$.
	
	\item $\gamma$ is a \emph{weakly-geodesic ray} if for every $y
	\in X$ and every $\varepsilon > 0$, there is an integer $N$
	such that
	\begin{equation*}
	  |d(\gamma(t), \gamma(0)) - t| < \varepsilon
	\end{equation*}
	and
	\begin{equation*}
	  |d(\gamma(t), y) - d(\gamma(s), y) - (t - s)| < \varepsilon
	\end{equation*}
	for all $t$, $s \in T$ with $t$, $s \ge N$.
    \end{enumerate}
\end{defn}

It is immediate that every geodesic ray is an almost-geodesic ray. 
Rieffel showed that every almost-geodesic ray is a weakly-geodesic
ray.  The significance of weakly geodesic rays is that they give the 
points on the metric boundary in reasonable metric spaces.

\begin{theorem}[Rieffel]\label{thm:metricboundary}
    Let $(X,d)$ be a complete, locally compact metric space, and let
    $\gamma: T \to X$ be a weakly geodesic ray in $X$.  Then
    \[
      \lim_{t \to \infty} f(\gamma(t))
    \]
    exists for every $f \in \mathcal{G}(X,d)$, and defines an element
    of $\partial_{d}X$.  Conversely, if $d$ is proper and if $(X,d)$
    has a countable base, then every point of $\partial_{d}X$ is
    determined as above by a weakly-geodesic ray.
\end{theorem}

This is similar in character to the definition of the Gromov boundary,
although the reliance on weakly-geodesic rays is necessary in general. 
Rieffel defined any point $\partial_{d}X$ which is the limit of an
almost-geodesic ray to be a \emph{Busemann point}, and it was shown
in~\cite{websterwinchester:busemann} that even for simple hyperbolic
spaces the metric boundary may have non-Busemann points.  It is an 
open question as to whether this phenomenon can occur with 
word-hyperbolic groups.

Unlike the Gromov boundary, the metric boundary is, in general,
dependent upon the choice of metric.  For example, different
generating sets for an infinite discrete group generally give distinct
metric boundaries for the corresponding word-length metrics.

From a practical viewpoint, the initial definition of the metric
boundary means that a sequence $x_{n} \in X$ converges to a point on
the metric boundary iff $x_{n}$ is eventually outside any compact
subset of $X$, and $\varphi_{y}(x_{n})$ converges for all $y \in X$. 
Two sequences converge to the same point on the metric boundary iff
\[
  \lim_{n \to \infty} \varphi_{z}(x_{n}) = \lim_{k \to \infty} 
  \varphi_{z}(y_{k})
\]
for every $z \in X$.  We can extend the functions $\varphi_{y}$ to the 
boundary by letting
\[
  \varphi_{y}(\omega) = \lim_{n \to \infty} \varphi_{y}(x_{n})
\]
for any sequence $x_{n} \to \omega \in \partial_{d}X$.  Then a
sequence $x_{n} \in \overline{X}^{d}$ converges to $x \in
\partial_{d}X$ iff $\varphi_{y}(x_{n}) \to \varphi_{y}(x)$ for all $y
\in X$, and this is sufficient to determine the topology of the metric
compactification.

\section{The Gromov Boundary as a Quotient}

We observe that the functions $\varphi_{y}$ and the inner product are
closely related, since
\[
  (x \cdot y) = \frac{1}{2}(\varphi_{y}(x) + d(y,0)),
\]
and furthermore, they play similar roles in the definitions of Gromov
and metric boundaries.  It is natural, therefore, to ask what
relationship there may be between the two different boundaries.

The key observation is that the triangle inequality implies that for
any $z \in X$,
\[
  (x \cdot y) \ge \frac{1}{2}(d(x,0) + d(y,0) - d(x,z) - d(y,z)) = 
  \frac{1}{2}(\varphi_{z}(x) + \varphi_{z}(y)),
\]
with equality iff $z$ lies on a geodesic path $[x,y]$.  We will want
to show that that $(x \cdot y)$ gets large for elements from various
sequences, and this implies that all we need do is find a $z$ so that
$\varphi_{z}(x) + \varphi_{z}(y)$ is large.

The following lemma tells us that as we get close to a metric boundary 
point, we can find $z$ such that $\varphi_{z}$ is large.

\begin{lemma}
    Let $(X,d)$ be a proper geodesic metric space with a distinguished 
    base-point 0.  Then for any $\omega$ in the metric boundary of 
    $X$, and any $N$, there is a point $z \in X$ such that 
    $\varphi_{z}(\omega) > N$.
\end{lemma}

\begin{proof}
    Let $x_{n}$ be any sequence which converges to $\omega$.
    
    Let $r > 0$ and consider a collection of minimal paths $[0,x_{n}]$
    for $n$ large enough that $d(0,x_{n}) > r$.  Because $(X,d)$ is a
    geodesic metric space, there must be a unique point $y_{n}$ in
    each of these paths in the sphere $S(0,r)$ of radius $r$, centred
    at 0.  Since $(X,d)$ is proper the sphere $S(0,r)$ is
    compact, and so given any $\varepsilon > 0$ we must be able to
    find at least one point $z_{r} \in S(0,r)$ such that an infinite
    number of the $y_{n}$ lie within $\varepsilon/2$ of $z_{r}$.  Let
    $x_{n_{j}}$ be the subsequence of $x_{n}$ corresponding to this
    infinite subset.  Then we have, for $r > \varepsilon$ and $j$
    sufficiently large,
    \[
      d(0,x_{n_{j}}) = d(0,y_{n_{j}}) + d(y_{n_{j}}, x_{n_{j}}) > 
      d(0,z_{r}) + d(z_{r},x_{n_{j}}) - \varepsilon,
    \]
    or, equivalently,
    \[
      \varphi_{z_{r}}(x_{n_{j}}) = d(0,x_{n_{j}}) - d(z_{r},x_{n_{j}})
      > d(0,z_{r}) - \varepsilon = r - \varepsilon.
    \]
    Taking limits, we conclude that
    \[
      \varphi_{z_{r}}(\omega) \ge r - \varepsilon.
    \]
    
    Hence, given any $N$, we can choose $r$ and $\varepsilon$ such 
    that $r - \varepsilon > N$, and obtain a point $z$ such that
    \[
      \varphi_{z}(\omega) > N.
    \]
\end{proof}

This lemma has two immediate corollaries:

\begin{corollary}
    Let $(X,d)$ be a proper geodesic metric space with a distinguished
    base-point $0$, and let $x_{n} \to \omega \in \partial_{d}X$. 
    Then $x_{n}$ converges to infinity in the Gromov sense.
\end{corollary}
\begin{proof}
    We know that for all $z$, $\varphi_{z}(x_{n})$ eventually gets
    close to $\varphi_{z}(\omega)$.  Hence by the previous lemma, for
    any $N$ can find a $z$ such that $\varphi_{z}(x_{n}) > N$ for all 
    $n$ sufficiently large.
    
    However, we than have that if $x_{n}$ and $x_{m}$ are large 
    enough that both $\varphi_{z}(x_{n})$ and $\varphi_{z}(x_{n})$ are 
    greater than $N$, then
    \[
      (x_{n} \cdot x_{m}) \ge \frac{1}{2}(\varphi_{z}(x_{n}) + 
      \varphi_{z}(x_{m})) > N
    \]
    
    Therefore
    \[
      \lim_{n,m \to \infty} (x_{n} \cdot x_{m}) = \infty
    \]
    and so $x_{n}$ goes to infinity in the Gromov sense.
\end{proof}

Let $(x_{n})$ and $(y_{k})$ be two sequences in $X$ which converge to
points on the metric boundary.  We will say that $(x_{n}) \sim_{d}
(y_{k})$ if these two sequences converge to the same metric boundary point. 
Similarly, if these sequences converge to points on the Gromov
boundary, we will say that $(x_{n}) \sim_{G} (y_{k})$.  Note that despite the
notation $\sim_{G}$ is not necessarily an equivalence relation.

\begin{corollary}
    Let $(X,d)$ be a proper geodesic metric space.  Then $(x_{n})
    \sim_{d} (y_{k})$ implies $(x_{n}) \sim_{G} (y_{k})$.
\end{corollary}

\begin{proof}
    Let $x_{n}$ and $y_{n}$ both converge to $\omega$.  Using the
    lemma, we can find a point $z$ so that $\varphi_{z}(\omega)$ is
    arbitrarily large, and since both $\varphi_{z}(x_{n})$ and
    $\varphi_{z}(y_{n})$ converge to $\varphi_{z}(\omega)$, for any
    number $N$ we can find $z$ such that both $\varphi_{z}(x_{n})$ and
    $\varphi_{z}(y_{n})$ are greater then $N$ for all $n$ sufficiently 
    large.
    
    Hence
    \[
      (x_{n} \cdot y_{n}) \ge \frac{1}{2}(\varphi_{z}(x_{n}) + 
      \varphi_{z}(y_{n})) > N
    \]
    for all $n$ sufficiently large, and so
    \[
      \lim_{n \to \infty} (x_{n} \cdot y_{n}) = \infty,
    \]
    and so $(x_{n}) \sim_{G} (y_{k})$.
\end{proof}

These two corollaries mean that we have a well-defined relation $\sim$
on $\partial_{d}X$ given by $\omega_{1} \sim \omega_{2}$ iff given any
$x_{n} \to \omega_{1}$ and $y_{k} \to \omega_{2}$, we have $(x_{n})
\sim_{G} (y_{k})$.  Furthermore, if $\sim_{G}$ is an equivalence
relation (as it is for hyperbolic spaces), then $\sim$ is an
equivalence relation on $\partial_{d}X$, and moreover $\partial_{G}X =
\partial_{d}X/\sim$ as sets.  As usual, we will denote the equivalence
class of a point $\omega$ in the metric boundary by $[\omega]$.

What we want is to show that we in fact have $\partial_{G}X =
\partial_{d}X/\sim$ as topological spaces.  In other words, we need 
to show that the quotient map is continuous.

\begin{lemma}
    Let $(X,d)$ be a proper hyperbolic metric space.  If $\omega_{n}
    \to \omega$ in $\partial_{d}X$, then $[\omega_{n}] \to [\omega]$
    in $\partial_{G}X$.
\end{lemma}
\begin{proof}
    Let $\delta > 0$ be the hyperbolic constant from
    (\ref{eqn:hyperboliccond}), $x_{k} \to \omega$ and $x_{n,k} \to
    \omega_{n}$.  We know that we can find $z$ such that
    $\varphi_{z}(\omega)$ is arbitrarily large, and since we have
    $\varphi_{z}(\omega_{n}) \to \varphi_{z}(\omega)$, we can choose
    $z$ such that $\varphi_{z}(\omega_{n})$ is also arbitrarily large,
    for all $n$ sufficiently large.  Indeed, as in the previous
    corollaries, we have, given any $N > 0$ we can find a number $M$
    such that for all $n > M$, there is a number $K_{n}$ such that
    \[
      (x_{n,k} \cdot x_{k}) \ge \frac{1}{2}(\varphi_{z}(x_{n,k}) +
      \varphi_{z}(x_{k})) > N + 2\delta
    \]
    for all $k > K_{n}$.
    
    Now if $y_{n,k} \to [\omega_{n}]$ and $y_{n} \to [\omega]$, we 
    know that we can find a subsequence of each sequence such that
    \[
      \liminf_{k \to \infty} (y_{n,k} \cdot y_{k}) = \lim_{k \to
      \infty} (y_{n,k_{j}} \cdot y_{k_{j}}).
    \]
    Furthermore, since $x_{n,j} \to [\omega_{n}]$ we conclude that for
    any $N$ there is some $J_{n}$ such that
    \[
      (y_{n,k_{j}} \cdot x_{n,j}) > N + 2\delta
    \]
    for all $j > J_{n}$, and similarly that there is some $J$ such that
    \[
      (y_{k_{j}} \cdot x_{j}) > N + \delta
    \]
    for all $j > J$.
    
    Hence, given any $N$, and fixing some $n > M$, then we have
    \[
      (y_{n,k_{j}} \cdot x_{j}) \ge \min \{(y_{n,k_{j}} \cdot 
      x_{n,j}), (x_{n,j} \cdot x_{j}) \} - \delta > N + \delta
    \]
    for all $j > \max \{J_{n}, K_{n}\}$.  But then
    \[
      (y_{n,k_{j}} \cdot y_{k_{j}}) \ge \min \{(y_{n,k_{j}} \cdot 
      x_{j}), (y_{k_{j}} \cdot x_{j}) \} - \delta > N
    \]
    for all $j > \max \{J_{n}, K_{n}, J\}$.  Therefore, for any $n > 
    M$,
    \[
      \lim_{k \to \infty} (y_{n,k_{j}} \cdot y_{k_{j}}) > N,
    \]
    and so
    \[
      \liminf_{k \to \infty} (y_{n,k} \cdot y_{k}) > N
    \]
    for all $n > M$.  And since $M$ does not depend on the choice of
    sequences converging to $[\omega_{n}]$ and $[\omega]$, we
    therefore have that
    \[
      ([\omega_{n}] \cdot [\omega]) = \inf\{\liminf_{k \to \infty}
      (y_{n,k} \cdot y_{k}) : y_{n,k} \to [\omega_{n}], y_{n} \to
      [\omega]\} > N
    \]
    for all $n > M$.
    
    Therefore
    \[
      \lim_{n \to \infty} ([\omega_{n}] \cdot [\omega]) = \infty,
    \]
    and so $[\omega_{n}] \to [\omega]$ in $\partial_{G}X$.
\end{proof}

So we have proved the following result.

\begin{theorem}\label{thm:quotientspace}
    Let $(X,d)$ be a proper, hyperbolic metric
    space.  Then there is a natural continuous quotient map from
    $\partial_{d}X$ onto $\partial_{G}X$.
\end{theorem}

\section{Boundaries of Word-Hyperbolic Groups}

We observe that if $G$ is a hyperbolic group, then the group acts on 
either boundary by taking a sequence $x_{k} \to \omega$ and letting
\[
  \alpha_{g}(\omega) = \lim_{k \to \infty} gx_{k}.
\]
This is a continuous action on either boundary.  Clearly the quotient
map is equivariant for these two actions, since if $\omega \sim
\omega'$, we can easily see that $\alpha_{g}(\omega) \sim
\alpha_{g}(\omega)$ by simply changing the base point of the inner
product to $g$.

An action of a topological group $G$ on a topological space $X$ is
\emph{amenable} if there is a net of continuous maps $(m_{\lambda}: 
X \to M_{1}^{+}(G))_{\lambda \in \Lambda}$, where $M_{1}^{+}(G)$ is 
the set of Borel probability measures on $G$, such that
\[
  \lim_{\lambda \in \Lambda} \| g \cdot m_{\lambda}(x) - 
  m_{\lambda}(g \cdot x) \| \to 0
\]
uniformly on compact subsets of $G \times X$.  Such a net of maps is
called an approximate invariant continuous mean.  It was shown by
E.~Germain (as discussed
in~\cite{anantharaman:amenablegroupoids,anantharaman:amenability})
that the action of a word-hyperbolic group $G$ on its Gromov boundary
is amenable.  Rieffel pointed out that if there were a continuous,
equivariant surjection from $\partial_{d}G$ to the Gromov boundary,
then the action of $G$ on the metric boundary must also be amenable. 
This is trivial given the above definition, since if $q: \partial_{d}G
\to \partial_{G}G$ is the quotient map of
Theorem~\ref{thm:quotientspace}, and $m_{\lambda}$ are the maps in an
approximate invariant continuous mean for the action of $G$ on
$\partial_{G}G$, then $m_{\lambda} \circ q$ are an approximate
invariant continuous mean for the action of $G$ on $\partial_{d}G$.

\begin{corollary}
    If $G$ is word-hyperbolic group with a finite generating set, 
    and $d$ is the word-length metric, then the group action on the 
    metric boundary is amenable.
\end{corollary}

This would seem to open the possibility of replicating Rieffel's work
on the metric boundary of $\integers^{d}$ in the setting of hyperbolic
groups.  However, Rieffel's procedure relied on the fact that the
action of $\integers^{d}$ on its metric boundary always has finite
orbits, and it seems unlikely that this criterion holds with any
frequency for general hyperbolic groups.

\bibliographystyle{hplain}

\bibliography{gromov}

\end{document}